# Max/Min Puzzles in Geometry II†

Revised

James M Parks
Dept. of Math.
SUNY Potsdam
Potsdam, NY
*parksjm@potsdam.edu*
January 20, 2022

**Abstract.** In this paper we continue the investigation of finding the max/min polygons which can be inscribed in a given triangle. We are concerned here with equilateral triangles. This may seem uninteresting or benign at first, but there are some surprises later.

*An inscribed polygon in a given triangle is a polygon which has all of its vertices on the triangle.* In a previous paper [3] we considered the inscribed polygons: parallelograms, rectangles, and squares. In this paper we investigate inscribed equilateral triangles in a given triangle. Figure 1 shows two examples of inscribed equilateral triangles in $\Delta ABC$.

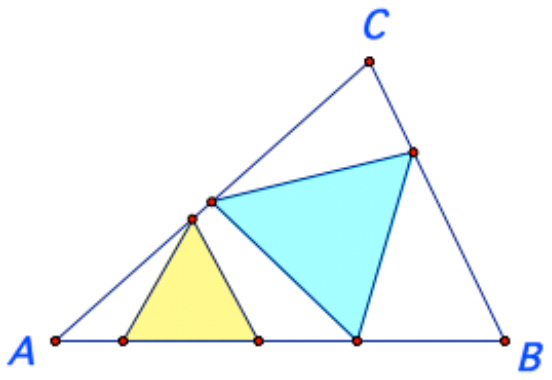

*Figure 1*

The first Puzzle is a quick investigation of properties of equilateral triangles (*ET* for short). All graphs were made with *Sketchpad**.

**Puzzle 1.** *If the sides of the ET $\Delta XYZ$ have known length s, determine the height h, and the area of $\Delta XYZ$ (denoted by (XYZ)). If the height h of ET $\Delta XYZ$ is known, determine the length of the side s, and the area (XYZ) of $\Delta XYZ$.*

Tip: The left and right hand halves of $\Delta XYZ$ are 30°-60°-90° triangles with sides in the ratios $s/2$, $s$, and $h = s\sqrt{3}/2$, Fig. 2.
Thus $(XYZ) = hs/2 = s^2\sqrt{3}/4$.
On the other hand, given $h$, the equality of ratios $h/\sqrt{3} = s/2$ holds.
Hence s $= 2h\sqrt{3}/3$, and $(XYZ) = h^2\sqrt{3}/3$.

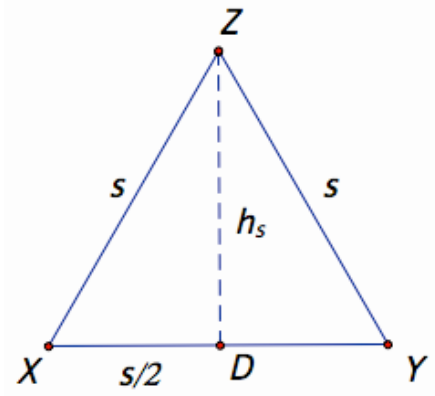

*Figure 2*

**Puzzle 2.** *Show that an inscribed ET of largest area in a given $\Delta ABC$, always shares at least one vertex with $\Delta ABC$.*

By trial and error it can be shown that there are only two ways an *ET* of maximum area will fit in a given $\Delta ABC$. Either with its base on one side of the triangle and its upper vertex coincides with the upper vertex of the given triangle, Fig. 3 Left, or it's base on one side of the triangle and at least one of it's base vertices coincides with one of the base vertices of the given triangle, Fig. 3 Right. That one side of the largest *ET* must be on a side of the triangle is well known [5], [6].

---

†Some of this material has appeared previously [4].
*All graphs in this paper were made with *Sketchpad v5.10 BETA.*

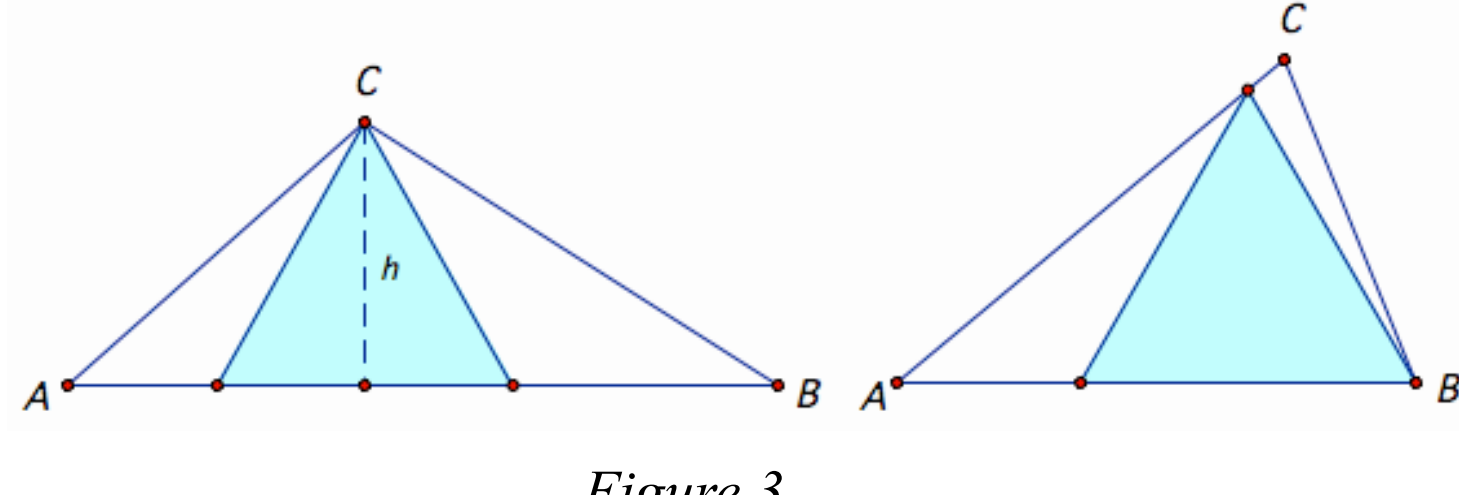

*Figure 3*

Observe that if the base angles of a side of $\Delta ABC$ are both greater than $60°$, then there is no inscribed *ET* possible on that side (look at examples, or see Fig. 13).
In Fig. 3 Left, notice that both base angles of $\Delta ABC$ are less than $60°$, so the area of the largest *ET* for this side is determined by the height h of $\Delta ABC$. From *Puzzle 1*, the area of the *ET* is $h^2\sqrt{3}/3$.
In Fig. 3 Right, base angle $B$ of $\Delta ABC$ is greater than $60°$, and base angle $A$ is less than $60°$, so it will take some additional calculating to compute the area of the largest *ET* for that side of $\Delta ABC$ (see *Puzzles 3 & 4*).

Since we are looking for the ***inscribed ET of maximum area***, if we can determine the largest inscribed *ET* for each side of the given triangle, then *the largest of these three ETs is called the* ***max ET*** *of the triangle*, and *the smallest of these 3 is called the* **min ET** *of the triangle*.
By the results in *Puzzle 2* we know exactly where to look for the *max* and *min ETs*. We will see below that this makes our search much easier.

We consider an example in *Puzzle 3* which illustrates Fig. 3 Right.
The angle notation is standard, angle CAB = $\alpha$, and so on.

**Puzzle 3**. *Let $\Delta ABC$ be the given triangle with $\alpha = 45°$, $\beta = 75°$, $\gamma = 60°$, and side a of length 2, Fig. 4. Let $\Delta BDE$ be the largest ET with base EB on side AB, and D on side AC. Compute s, the length of the side of $\Delta BDE$, and the exact area of $\Delta BDE$. Is $\Delta BD$ the max ET in $\Delta ABC$?*

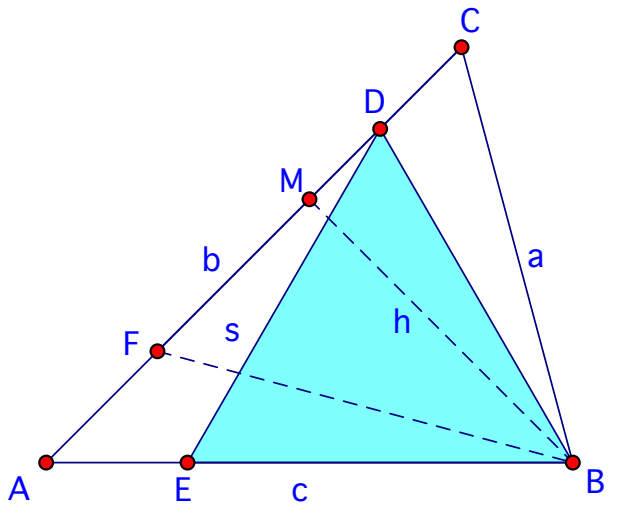

*Figure 4*

Hint: Drop a perpendicular from $B$ to $AC$ at $M$.
Then $\Delta MBC$ is a $30°$-$60°$-$90°$ triangle with $|BM| = \sqrt{3}$.
This means that $c = \sqrt{6}$, since $\Delta ABM$ is a $45°$-$45°$-$90°$ triangle.
By the Law of Sines, $s = \sqrt{6}(\sin 45°/\sin 75°) \sim 1.793$, since angle $BDA = 75°$. Hence $(BDE) = s^2\sqrt{3}/4 \sim 1.392$.
For side $b$ construct a point $F$ on $AC$ such that $|FC| = 2$.
Then $\Delta BCF$ is an *ET* which is the largest *ET* on both sides $a$ and $b$, and $(BCF) = \sqrt{3} \sim 1.732$, by *Puzzle 1*.
Thus $\Delta BCF$ is the *max ET* in $\Delta ABC$.

The key to solving *Puzzle 3* was knowing that $\gamma = 60°$. But what if $\gamma \neq 60°$? How do we compute the area of $\Delta BDE$ in that case?
Here's a way to compute the area of such an *ET*.

**Puzzle 4.** *Compute the area of the largest inscribed ET on side c of* $\Delta ABC$, *for* $\alpha < 60°$, $\beta > 60°$.

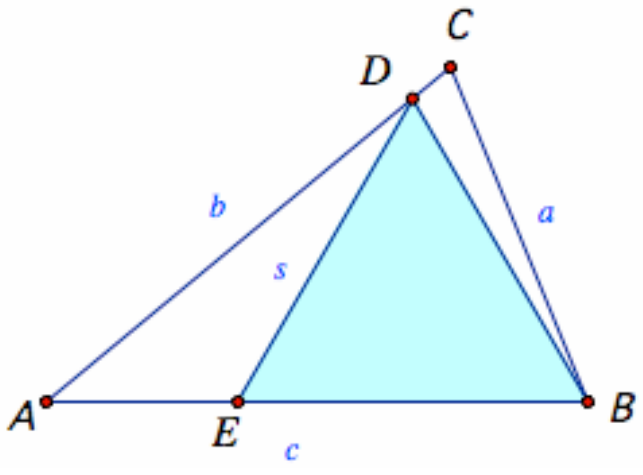

*Figure 5*

Hint: By the Exterior Angle Theorem, angle $CDB = 60°+\alpha$.
By the Law of Sines, $s = a\sin\gamma/\sin(60°+\alpha)$.
But if $h_b$ is the altitude of $B$ to $AC$ at $M$, $h_b = BM$, as in Fig.4, then $\sin\gamma = h_b/a$, hence $s = h_b/\sin(60° + \alpha)$.
So, by the solution to *Puzzle 1*, $(BDE) = (h_b/\sin(60° + \alpha))^2\sqrt{3}/4$.

*We will now determine the location of the max and min ETs for all triangles.*
Also, ***we assume from here on that*** $\Delta ABC$ ***is scalene, unless noted otherwise****, and recall in particular that* $\alpha > \beta > \gamma$ *iff* $a > b > c$, [2, Vol.1, Props. 18, 19].

By symmetry, any result we get for $\Delta ABC$ under these assumptions will also hold for a triangle reflected about the vertical line through the mid-point of the base $AB$. So we need only consider those triangles with upper vertex $C$ to the left of this vertical line. If $C$ is on this line, then $\Delta ABC$ is an isosceles triangle, which is a simpler case to work with.

In the chart in Fig. 6 there is a list of examples generated on *Sketchpad*. The ratio of the area of the *ET* to the area of $\Delta ABC$ is listed in the right-hand column, marked *Ratios*. This makes it easier to compare the sizes of the areas of the *ETs*. The bold **x** marks the side *a, b, c,* where the *max ET* occurs, and the italic *x* marks the side where the *min ET* is found. Several patterns can be observed here.

The general results are given in *Puzzle 5*. In each part of *Puzzle 5* an example(s) of triangles from this chart is considered.

| Angles<br>α - β - γ | Sides: long/middle/short<br>a | b | c | Ratios |
|---|---|---|---|---|
| 70-65-45 | | *x* | | *0.638* |
| " | **x** | | | **0.677** |
| 80-65-35 | | *x* | | *0.461* |
| " | **x** | | | **0.547** |
| 90-70-20 | | *x* | | *0.288* |
| " | **x** | | | **0.327** |
| 90-80-10 | | *x* | | *0.168* |
| " | **x** | | | **0.173** |
| 75-60-45 | | *x* | | *0.589* |
| " | **x** | | | **0.732** |
| 90-60-30 | | *x* | | *0.378* |
| " | **x** | | **x** | **0.500** |
| 100-60-20 | | *x* | | *0.265* |
| " | **x** | | **x** | **0.348** |
| 90-50-40 | | *x* | | *0.437* |
| " | **x** | | | **0.569** |
| 110-40-30 | | *x* | | *0.297* |
| " | **x** | | | **0.395** |
| 120-35-25 | | *x* | *x* | *0.243* |
| " | **x** | | | **0.325** |
| 130-30-20 | | | *x* | *0.193* |
| " | **x** | | | **0.258** |
| 140-25-15 | | | *x* | *0.150* |
| " | **x** | | | **0.193** |

*Figure 6*

**Puzzle 5.** *Given* $\Delta ABC$, *suppose the sides a, b, c satisfy* $a > b > c$. *Show that the max ET is always on the long side a (a and c, when* $\beta = 60°$*), and the min ET is on:*

***A****. side c, when* $120° < \alpha$*; and sides b and c, when* $\alpha = 120°$*;*
***B****. side b, when* $60° < \alpha < 120°$*, and* $\beta < 60°$*;*
***C****. side b, when* $60° \leq \beta$.

In Figures *7 - 10*, notice that all *3* largest *ETs* (*2 ETs* coincide in Fig.10) have vertex *A* in common. Then *ET* $\Delta GHA$, on the long side *a*, is the *max ET*, since *JA* and *FA* are interior to it (a consequence of Euclid Prop. 19 [2]). Since $\Delta GHA$ and $\Delta AEF$ coincide when $\beta = 60°$, Fig. 10, the *max ET* is on both the long side *a* and the short side *c* for this case. If $\beta > 60°$, Fig. 10, then $\alpha > \beta > 60°$, and $m\angle JAB > m\angle IBA$. So $|HI| > |JA|$, by Euclid Prop. 19 [2], and $\Delta GHI$ is the *max ET* on the long side *a.*
These results agree with what is observed in the chart in Fig. 6.

**A.** Fig. 7 is a *130°-30°-20°* triangle with *ET* $\Delta GHA$ on the long side *a*, which is clearly the *max ET,* since $|AG| > |AJ|$ and $|AG| > |AF|$. If we compare $|AJ|$ and $|AF|$ we discover that $\Delta AEF$ is the *min ET*.
For the general case for $\alpha > 120°$, since $m\angle CAJ = 60° = m\angle FAB$, it follows that $m\angle AJC > m\angle AFB$. So, in $\Delta AJF$, $s_b = |JA| > |FA| = s_c$, and the *min ET* is on the short side *c*.

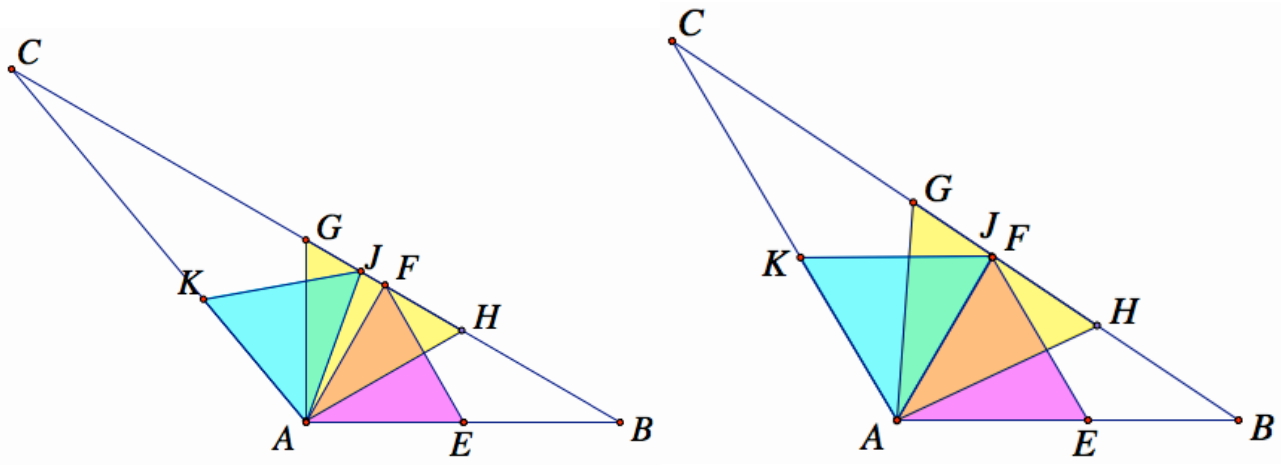

*Figure 7* *Figure 8*

Figure 8 is a *120°-35°-25°* triangle with the *max ET*, $\Delta GHA$, on the long side *a*. From the first of part **A** it follows that $\Delta GHA$ is larger than the other two *ETs*, and $\Delta AEF$ and $\Delta JKA$ are congruent. For a general argument, if $\alpha = 120°$, then $s_b = |JA| = |FA| = s_c$, so $\Delta AEF$ and $\Delta AJK$ are congruent, and the *min ET* is on both the short side *c,* and the middle side *b*.

**B.** In the *110°-40°-30°* triangle, Fig. 9, the area of the *min ET,* $\Delta AJK$ on side *b*, is less in area than *ET* $\Delta AEF$ on side *c*, which is less in area than the *max ET* $\Delta GHA$ on side *a,* by comparing sides *AG, AF,* and *AJ*.
For a formal argument, the situation is similar to part **A,** since $m\angle CAJ = 60° = m\angle FAB$, so we have $m\angle AJC > m\angle AFB$. In $\Delta AJF$ it follows that $s_b = |JA| < |FA| = s_c$, and $\Delta AJK$ is the *min ET* on the middle side *b*.

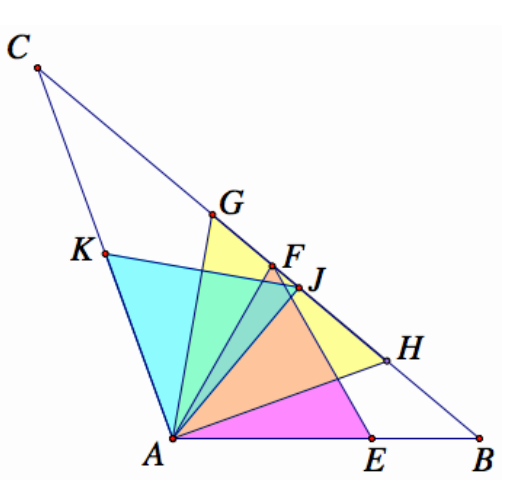

*Figure 9*

**C.** In the *75°-60°-45°* triangle Fig. 10, the *ETs* on sides *c* and *a* coincide, and have sides of length *c*, while the *ET* on side *b,* $\Delta AJK$, has side $s_b < c$.
In general, if $\beta = 60°$, then $\Delta AEF$ and $\Delta GHA$ coincide. Also *AJ* is interior to $\Delta AEF$, $s_b = |JA| < |AE| = s_c = c$. This means the *min ET* is $\Delta AJK$ on the middle side *b*.

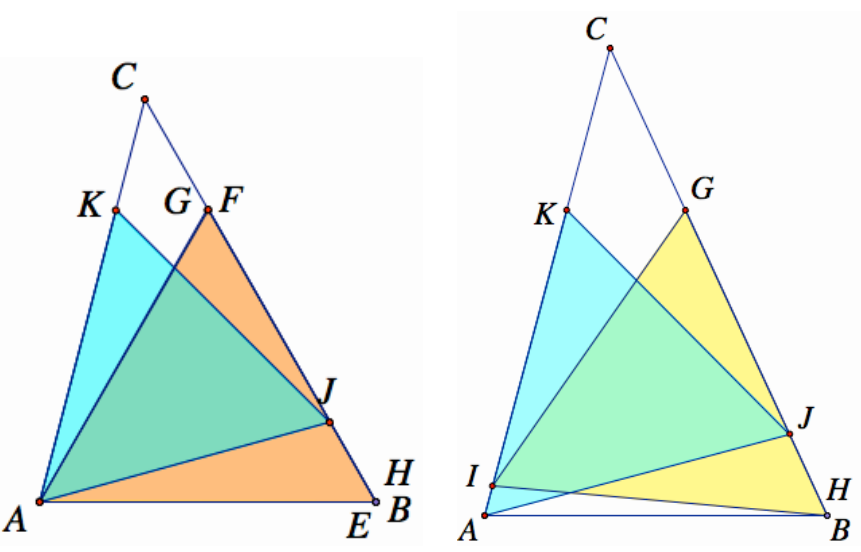

*Figure 10* *Figure 11*

In the *70°-65°-45°* triangle, Fig.11, the length of the side of the *min ET*, $\Delta AJK$ on side *b*, is less than that of the *max ET,* $\Delta GHI$ on side *a,* by comparing $|AJ|$ and $|IB|$.
For the general argument, suppose $\beta > 60°$, then $\alpha > \beta > 60°$, and by subtracting *60°* from $\alpha$ and $\beta$ it follows that $m\angle JAB > m\angle IBA$. This means $s_a = |HI| > |JA| = s_b$, and $\Delta JKA$ is the *min ET* on the middle side $b$. There is no *ET* on side $c$ in this case, as the base angles are greater than *60°*.

In Fig. 12 the graphs of the regions which correspond to Puzzles ***5A***, ***5B***, and ***5C***, are shown. The regions are labeled with the puzzle letter. A typical triangle $\Delta ABC$ is shown with vertex $C$ in the region ***B***. Ignoring isosceles triangles for now, the dashed lines in red and green are off limits, as is the white circular area (recall it is assumed that $a > b > c$).

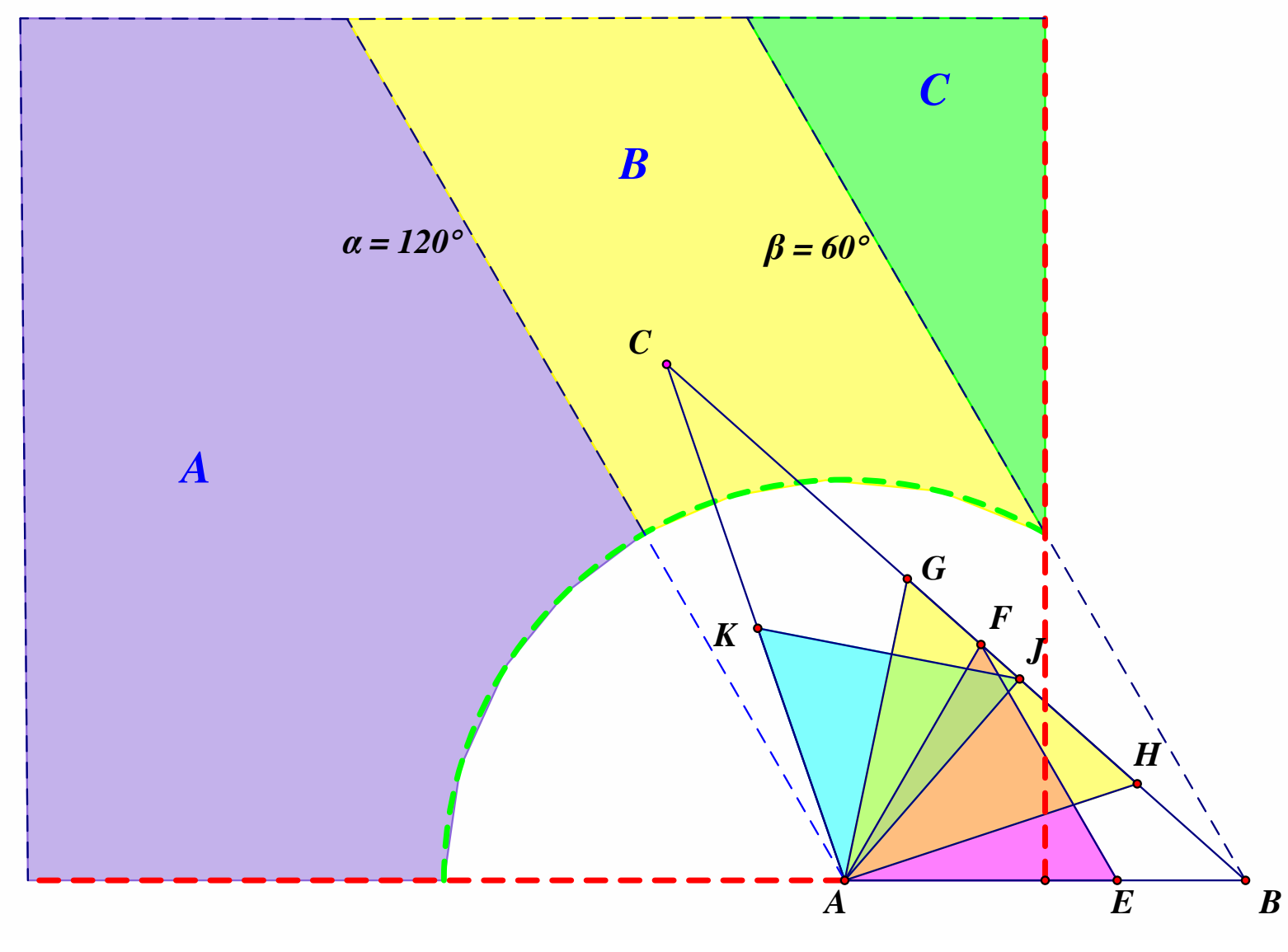

*Figure 12*

If *ETs* of largest area on side $c$ are allowed, when vertex $C$ is in region ***C***, there is a *wedged ET* (*WET* for short) [1], [4].
These are not inscribed *ETs* as defined above, since the *ETs* do not have all *3* vertices on $\Delta ABC$, but they are the largest *ETs* on side $c$. The term *wedged* was first used by Calabi in [1].

The new results under this change will happen in *Puzzle* ***5C***, for then a *WET* on the base side $c$ of the triangle can be included, see Figures 11 and 13.
We list these new results to *Puzzle 5C* here as *5C.1, 5C.2,* and *5C.3*.

**5C.1.** *Let $\Delta ABC$ satisfy $60° \leq \beta < 80°$, $60° < \alpha < 90°$, and $\alpha + \beta/2 < 120°$. Then the max WET is on the short side c, and the min WET is on the middle side b.*

**5C.2.** *In $\Delta ABC$, if we assume:*

*i.) $60° < \beta < 80° < \alpha < 90°$, $\alpha/2 + \beta < 120°$, and $\alpha + \beta/2 \geq 120°$, then the min WET is on the middle side b, and the max WET is on both the short side c and the long side a, when $\alpha + \beta/2 = 120°$, and on the long side a, when $\alpha + \beta/2 > 120°$.*

*ii.) $75° < \beta < 90°$, $80° < \alpha < 90°$, and $\alpha/2 + \beta \geq 120°$, then the max WET is on the long side a, and the min WET is on the short side c, when $\alpha/2 + \beta > 120°$; and on the short side c and the middle side b, when $\alpha/2 + \beta = 120°$.*

**5C.3.** *In $\Delta ABC$, if $60° < \beta < 90° < \alpha$, then the max WET is on the long side a, and the min WET*

*i.) is on the middle side b, when $\alpha/2 + \beta < 120°$, and on the short side c and the middle side b, when $\alpha/2 + \beta = 120°$,*

*ii.) is on the short side c, when $\alpha/2 + \beta > 120°$.*

As might be expected, the arguments for these results are a little more involved than those given for Puzzles ***5A**, **B**,* and ***C*** above. For example, the following is an argument for result ***5C.2.i***.

Suppose $60° < \beta < 80° < \alpha < 90°$, $\alpha/2 + \beta < 120°$, $\alpha + \beta/2 = 120°$, and $20° < \gamma < 30°$, by the boundaries of the subregion. This means $90° < m\angle AEC$, and $m\angle BGC < 100°$. Since $\beta < \alpha$, $m\angle EAB > m\angle GBA$, and $|GB| > |AE|$. So $\Delta BHG$ is a possible *max WET* on side $a$, and $\Delta AEF$ is a possible *min WET* on side $b$.

*Figure 13*

Let $BJ$ be the angle bisector of $\beta$, $J$ on $CA$, Fig. 13. Then $m\angle AJB = 60°$, since $\alpha + \beta/2 = 120°$. This also means $\gamma + \beta/2 = 60°$, by the Ext. Angle Th., so $m\angle CJB = 120° = m\angle CHG$, and $m\angle CGH = \beta/2$. But then $m\angle AGB = \alpha$, since $\alpha + \beta/2 = 120°$, so $\Delta ABG$ is isosceles, and $c = |GB|$. Therefore $\Delta ABD$ and $\Delta BHG$ are equal *max WETs* on sides $c$ and $a$, resp., and $\Delta AEF$ is the *min WET* on side $b$.

On the other hand, if $60° < \beta < 80° < \alpha < 90°$, $\alpha/2 + \beta < 120°$, and $\alpha + \beta/2 > 120°$, then $15° < \gamma < 30°$ (check the conditions for $\alpha = 90°$), and as above $\Delta CAE$ and $\Delta CBG$ are similar, so $m\angle AGB = m\angle AEB$. Then $\alpha > \beta$ implies $\alpha + \beta/2 > \alpha/2 + \beta$, $m\angle EAB > m\angle GBA$, and $m\angle AEB > \beta$, so $|GB| > c > |AE|$. This means that $\Delta BHG$ is the *max WET* on side $a$ and $\Delta AEF$ is the *min WET* on side $b$.

The argument for ***5C.2.ii*** is analogous, but the angle $\alpha$ is bisected instead of angle $\beta$ in that case. Also, a similar argument can be used in ***5C.3***.

The other cases for isosceles triangles with central angle less than $90°$ are contained in results ***5C.1***, and ***5C.2.ii***.

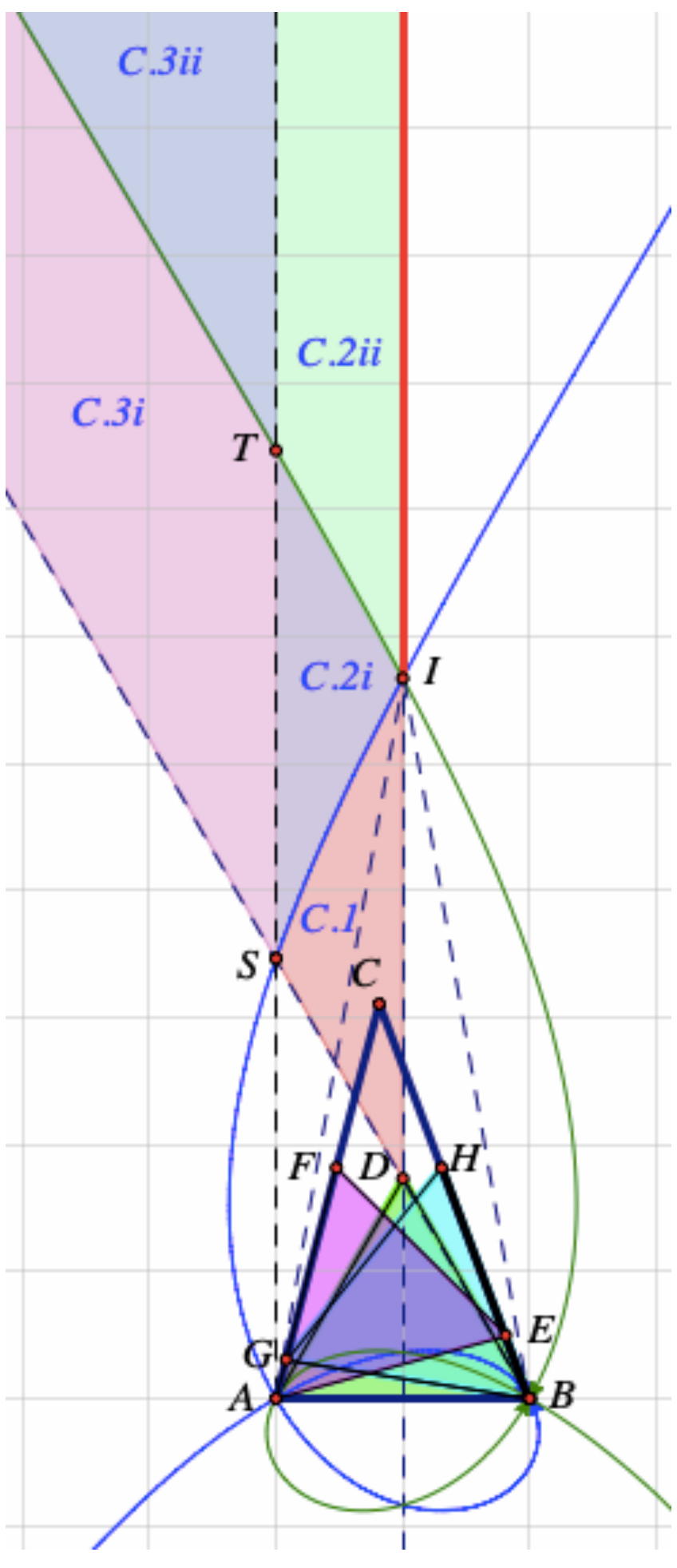

*Figure 14*

These new results for Puzzle ***5C*** divide region ***C*,** Fig. 12, into subregions which correspond to each of these new results, Fig. 14.

The subregion ***C.2.i*** is of particular interest, as two of the boundaries are the graphs of cubic equations, determined by the conditions on the angles in that part.

In region ***C.2.ii***, if $\alpha/2 + \beta = 120°$, then $\gamma = 60° - \alpha/2$, so by the law of sines, if $c = 1$, then $b = \sin\beta/\sin\gamma =$ $\sin(120° - \alpha/2)/\sin(60° - \alpha/2) = \sin(60° + \alpha/2)/\sin(60° - \alpha/2)$, [2].

Converting to polar coordinates: $r = \sin(60° + \theta/2)/\sin(60° - \theta/2)$, for $80° < \theta < 90°$, $A = (0, 0)$ the polar pole, $B = (1, 0)$, $\theta = \alpha$, and $AB$ on the polar axis.

Changing to rectangular coordinates, by the usual method, determines the associated cubic equation:

$$\sqrt{3}x^3 + x^2(y - 2\sqrt{3}) + \sqrt{3}x(y^2 + 1) + y^3 - 2\sqrt{3}y^2 - y = 0.$$

The graph of this curve (green) is the lower bound of region ***C.2.ii*** (green) and the upper bound of region ***C.2.i***.

Similarly, the lower bound of region ***C.2.i*** (blue) which is the upper bound of region ***C.1*** satisfies the condition $\alpha + \beta/2 = 120°$, so $\gamma = 60° - \beta/2$.

By the law of sines, if $c = 1$, then $b = \sin\beta/\sin\gamma = \sin(240° - 2\alpha)/\sin(\alpha - 60°)$.

The corresponding polar coordinates equation is: $r = \sin(240° - 2\theta)/\sin(\theta - 60°)$, for $80° < \theta < 90°$.

In rectangular coordinates this corresponds to the cubic equation:

$$\sqrt{3}x^3 - x^2(y + \sqrt{3}) + x(\sqrt{3}y^2 + 2y) + \sqrt{3}y^2 - y^3 = 0.$$

The graph of this curve (blue) is the lower bound of region ***C.2.i***, and the upper bound of ***C.1***.

The two cubic curves are symmetric about the vertical line through points *D* and *I*. The point labeled *I* is a point where the *2* cubics intersect, and also the upper vertex of a triangle which is a *WET* analog of Calabi's result for wedged squares in triangles, *Puzzle 6,* [1] [3].

It should be obvious that if $\Delta ABC$ is itself an *ET*, then all inscribed *ETs* of largest area will coincide with $\Delta ABC$. From our analysis above this is the only triangle where the 3 largest inscribed *ETs* are equal. However, if *WETs* are admitted, then a new possibility exists.

**Puzzle 6.** *If $\Delta ABC$ is isosceles with base angles of 80°, show that the largest WETs on sides a, b, and c are congruent.*

Hint: Let $\Delta BHG$, $\Delta AEF$, and $\Delta ABD$ be the largest *WETs* on sides *a, b,* and *c,* respectively, Fig. 15.

Since $\alpha = \beta = 80°$, it follows that $\gamma = 20°$, and also $m\angle EAB = m\angle GBA = 20°$.

So, $m\angle AEB = m\angle BGA = 80°$.

Thus, it follows that triangles $\Delta ABE$ and $\Delta ABG$ are isosceles, so $|GB| = |AB| = |AE|$.

Hence, all *3* largest *WETs:* $\Delta ABD$, $\Delta AEF$, and $\Delta BHG$, are congruent.

*Figure 15*